\pdfoutput=1
\RequirePackage{ifpdf}
\ifpdf 
\documentclass[pdftex]{sigma}
\else
\documentclass{sigma}
\fi

\usepackage{ifthen}

\newif\ifdraftversion
\draftversionfalse

\newif\ifarxivversion
\arxivversiontrue

\newif\iffinal
\finaltrue

\iffinal
\draftversionfalse
\fi

\usepackage{enumerate}
\usepackage{mathtools}

\makeatletter
\newsavebox{\@brx}
\newcommand{\llangle}[1][]{\savebox{\@brx}{\(\m@th{#1\langle}\)}%
 \mathopen{\copy\@brx\mkern2mu\kern-0.9\wd\@brx\usebox{\@brx}}}
\newcommand{\rrangle}[1][]{\savebox{\@brx}{\(\m@th{#1\rangle}\)}%
 \mathclose{\copy\@brx\mkern2mu\kern-0.9\wd\@brx\usebox{\@brx}}}
\makeatother

\usepackage{slashed}

\usepackage{cancel}

\newcommand{\vertiefetextbox}[2]{\settoheight{\meineboxheight}{#2}\addtolength{\meineboxheight}{#1}\lower#1\vbox to \meineboxheight{\hbox{#2}\vfil}}

\ifdraftversion
\newcounter{mnotecount}[section]
\renewcommand{\themnotecount}{\thesection.\arabic{mnotecount}}
\newcommand{\mnote}[1]{\protect{\stepcounter{mnotecount}}${\raisebox{0.5\baselineskip}[0pt]{\makebox[0pt][c]{\tiny\em{\red{$\bullet$\themnotecount}}}}}$\marginpar{\raggedright\tiny\em $\!\!\!\!\!\!\,\bullet$\themnotecount: #1}\ignorespaces}
\else
\newcommand{\mnote}[1]{}
\fi

\ifdraftversion
\newcommand{\mbernd}[1]{\mnote{B: #1}}
\newcommand{\mmattias}[1]{\mnote{M: #1}}
\newcommand{\msamuel}[1]{\mnote{S: #1}}
\newcommand{\fbernd}[1]{\footnote{B: #1}}
\newcommand{\fmattias}[1]{\footnote{M: #1}}
\newcommand{\fsamuel}[1]{\footnote{S: #1}}
\newcommand{\bernd}[1]{\red{B: #1}}
\else
\newcommand{\mbernd}[1]{}
\newcommand{\mmattias}[1]{}
\newcommand{\msamuel}[1]{}
\newcommand{\fbernd}[1]{}
\newcommand{\fmattias}[1]{}
\newcommand{\fsamuel}[1]{}
\newcommand{\bernd}[1]{}
\fi

\usepackage{xcolor}

\definecolor{darkgreen}{rgb}{0,0.6,0}
\definecolor{darkred}{rgb}{0.7,0,0}
\definecolor{darkblue}{rgb}{0,.1,.6}
\definecolor{darksharpblue}{rgb}{.4,.2,.8}
\definecolor{darkgray}{rgb}{0.3,.3,.3}
\definecolor{lightgreen}{rgb}{.9,1,.9}

\newcommand\red[1]{\textcolor{red}{#1}}

\let\O\relax
\let\div\relax
 \let\Re\relax

\let\mod\undefined

\let\epsilon\varepsilon

\let\theta\vartheta

\def\phi{{\varphi}}

\let\na\nabla
\newcommand\ceq{\coloneqq}

\let\witi\widetilde

\newcommand\grad{\operatorname{grad}}

\DeclareMathOperator{\tr}{tr}
\DeclareMathOperator{\id}{id}
\DeclareMathOperator{\spec}{spec}
\DeclareMathOperator{\univ}{univ}
\newcommand\Cl{\mathop{\mathrm{Cl}}}

\newcommand\rank{\mathop{\mathrm{rank}}}

\newcommand\image{\mathop{\mathrm{image}}}
\newcommand\upd{\mathrm{d}}

\newcommand\ball{ball}
\newcommand\Dirac{\slashed{D}} 

\newcommand\Dminimal{$\Dirac$-minimal}
\newcommand\enerimp{Q}

\newcommand\rangleH{\rangle_\HH}
\newcommand\rangleC{\rangle_\CC}
\newcommand\rangleR{\rangle_\RR}
\newcommand\rrangleH{\rrangle_\HH}

\newcommand\<{\langle}
\renewcommand\>{\rangle}

\newlength{\meineboxheight}
\newlength{\meineboxdepth}

\newlength{\templength}
\newcommand{\setzeboxnachlinksmath}[1]{\settowidth{\templength}{$#1$}\kern-\templength#1}
\newlength{\qedskiplength}
\newcommand\point{\mathrm{pt}}

\newcommand\datver[1]{\def\datverp
{\par\boxed{\boxed{\text{Version: #1; Run: \today}}}}}\datver{0.1}

\newcommand{\CC}{\mathbb C}
\newcommand{\HH}{\mathbb{H}}
\newcommand{\KK}{\mathbb K}

\newcommand{\NN}{\mathbb N}

\newcommand{\RR}{\mathbb R}

\newcommand{\ZZ}{\mathbb Z}

\let\Z\ZZ
\let\C\CC

\newcommand{\End}{\operatorname{End}}
\newcommand{\GL}{{\operatorname{GL}}}
\newcommand{\O}{{\operatorname{O}}}
\newcommand{\SO}{{\operatorname{SO}}}
\newcommand{\Spin}{{\operatorname{Spin}}}

\newcommand{\KO}{{\operatorname{KO}}}

\newcommand{\dvol}{\operatorname{dvol}}
\newcommand{\dV}{\,\mathrm{d}V}
\newcommand{\codim}{\operatorname{codim}}
\newcommand{\div}{\operatorname{div}}
\let\ol\overline

\newcommand\Sym[1][2]{\mathop{\mathrm{Sym}}\nolimits^{#1}}
\DeclareMathOperator{\spann}{span}
\DeclareMathOperator{\Re}{Re}

\newcommand{\mod}{\,\,\mathrm{mod}\,\,}

\newcommand{\maA}{\mathcal A}

\newcommand{\maE}{\mathcal E}

\newcommand{\maI}{\mathcal I}

\newcommand{\maM}{\mathcal M}
\newcommand{\maN}{\mathcal N}

\newcommand{\maR}{\mathcal R}
\newcommand{\maS}{\mathcal S}

\newcommand{\maU}{\mathcal U}
\newcommand{\maV}{\mathcal V}
\newcommand{\maW}{\mathcal W}

\newcommand\bullette{{\scalebox{0.5}{$\bullet$}}}
\newcommand\bullettette{{\scalebox{0.3}{$\bullet$}}}
\newcommand\cliff{\,\raisebox{.18\keyheightlength}{\bullette}\,}
\newcommand\argu{\,\raisebox{.1\keyheightlength}{\bullettette}\,}

\let\biargu\arrgu

\newtheorem{Theorem}{Theorem}[section]
\newtheorem*{Theorem*}{Theorem}

\newtheorem{Lemma}[Theorem]{Lemma}
\newtheorem{Proposition}[Theorem]{Proposition}

\theoremstyle{definition}

\newtheorem{Example}[Theorem]{Example}
\newtheorem{Remark}[Theorem]{Remark}
\newtheorem{Remarks}[Theorem]{Remarks}

\newlength{\keyheightlength}

\let\define\emph

\newcommand{\meineurl}[1]{\href{#1}{(URL)}}

\newcommand\eg{e.g.,\ \ignorespaces}
\newcommand\ie{i.e.,\ \ignorespaces}

\newcounter{proofstep}

\makeatletter
\newcommand{\xdashrightarrow}[2][]{\ext@arrow 0359\rightarrowfill@@{#1}{#2}}
\def\rightarrowfill@@{\arrowfill@@\relax\relbar\rightarrow}
\def\arrowfill@@#1#2#3#4{%
 $\m@th\thickmuskip0mu\medmuskip\thickmuskip\thinmuskip\thickmuskip
 \relax#4#1
 \xleaders\hbox{$#4#2$}\hfill
 #3$%
}
\makeatother
\newcommand{\dashto}{\xdashrightarrow{\kern10mm}}

\let\tab\=
\let\=\undefined
\newcommand\={\,=\,}

\newcommand\Rmin{\maR_{\mathrm{min}}}
\newcommand\Rcf{\maR_{\mathrm{cf}}}
\newcommand\Rinv{\maR_{\mathrm{inv}}}
\newcommand\Rpsc{\maR_{\mathrm{psc}}}

\numberwithin{equation}{section}

\begin{document}

\allowdisplaybreaks

\newcommand{\arXivNumber}{2508.01420}

\renewcommand{\PaperNumber}{102}

\FirstPageHeading

\ShortArticleName{The Space of Dirac-Minimal Metrics is Connected in Dimensions~2 and~4}

\ArticleName{The Space of Dirac-Minimal Metrics\\ is Connected in Dimensions~2 and~4}

\Author{Bernd AMMANN~$^{\rm a}$ and Mattias DAHL~$^{\rm b}$}

\AuthorNameForHeading{B.~Ammann and M.~Dahl}

\Address{$^{\rm a)}$~Fakult\"at f\"ur Mathematik, Universit\"at Regensburg, 93040 Regensburg, Germany}
\EmailD{\mail{bernd.ammann@mathematik.uni-regensburg.de}}
\URLaddressD{\url{https://www.berndammann.de}}

\Address{$^{\rm b)}$~Institutionen f\"or Matematik, Kungliga Tekniska H\"ogskolan, 100 44 Stockholm, Sweden}
\EmailD{\mail{dahl@kth.se}}
\URLaddressD{\url{https://www.kth.se/profile/dahl}}

\ArticleDates{Received August 11, 2025, in final form November 25, 2025; Published online December 06, 2025}

\Abstract{Let $M$ be a closed connected spin manifold. Index theory provides a topological lower bound on the dimension of the kernel of the Dirac operator which depends on the choice of Riemannian metric. Riemannian metrics for which this bound is attained are called Dirac-minimal. We show that the space of Dirac-minimal metrics on $M$ is connected if $M$ is of dimension~2 or~4.}


\Classification{53C27; 19K56; 58C40; 58J5}

\section{Introduction}

On a closed connected spin manifold $M$, index theory for the Dirac operator, in particular the family Atiyah--Singer index theorem is a strong and interesting tool to determine non-trivial topology of the space $\Rpsc(M)$ of Riemannian metrics with positive scalar curvature (psc metrics). However, when the index theoretical information vanishes, it remains a notoriously difficult problem to prove connectedness properties for $\Rpsc(M)$. For $\dim M=3$, a recent result by Bamler and Kleiner~\cite{bamler.kleiner:2019p} building on previous work by Marques~\cite{marques:2012} shows that $\Rpsc(M)$ is either empty or contractible.

The above mentioned applications of index theory rely on the fact that the Dirac operator for a psc metric is an invertible operator. Thus, the space $\Rinv(M)$ of metrics whose Dirac operator is invertible contains $\Rpsc(M)$. As a consequence, non-trivial homotopy groups in~$\Rpsc(M)$, that are detected by (untwisted) index theoretical methods remain non-trivial in the space~$\Rinv(M)$. For example, if the index $\alpha(M)\in \KO_m(\point)$, $m=\dim M$, does not vanish, we have~${\Rpsc(M)\subset\Rinv(M)=\varnothing}$. Similarly, spectral flow of the Dirac operator, or equivalently a family index theorem, can be used to show that $\Rpsc(M)$ and $\Rinv(M)$ have infinitely many components if $m\equiv 3 \mod 4$, $m\geq 7$, and at least two components, if $m\equiv 0,1\mod 8$. Such results use the fact that $\KO_{m+1}(\point)$ is non-trivial in these dimensions~$m$.

However, if the $\KO_{m+1}(\point)$-valued spectral flow between Riemannian metrics $g_1,g_2\!\in\! \Rpsc(M)$ is zero, it is difficult to show that $g_1$ and $g_2$ are in the same connected component of $\Rpsc(M)$. The main result of this article is a proof of the connectedness of $\Rinv(M)$ in the cases $m=2$ and $m=4$. Note that this is the first such connectedness result for $\Rinv(M)$. Our result is in contrast to the setting of metrics with positive scalar curvature in the case $m=4$: Seiberg--Witten theory tells us that on many connected closed oriented $4$-manifolds $\Rpsc(M)$ is not connected \cite[Sections~5 and~6]{ruberman:2001}.

If the index $\alpha(M)$ does not vanish, the space $\Rinv(M)$ is empty, and our setting has to be adapted. Let us recall the Atiyah--Singer index theorem for a closed spin manifold $M$ with its extension using Hitchin's $\alpha$-invariant in $\KO_m(\point)$. Let $\maR(M)$ be the space of Riemannian metrics on $M$.
For any $g\in \maR(M)$ we consider the spinor bundle $\Sigma M=M\times_\rho\Sigma_m$ that arises by considering a real irreducible representations $\rho\colon\Cl_m\to \End(\Sigma_m)$ of the $m$-dimensional real Clifford algebra, $m=\dim M$. Let us summarize some classical facts about such representations, as explained and proven, \eg in \cite[Theorem~I.5.8]{lawson.michelsohn:89}. In dimensions $m\not\equiv 3 \mod 4$, the choice of irreducible representation is unique up to isomorphism (in the sense of real representations).
The representation $\Sigma_m$ carries a complex structure, commuting with the $\Cl_m$-action if ${m\equiv 1,2,3,4,5\mod 8}$, and this complex structure is unique (up to isomorphisms in the sense of complex representations) for $m\equiv 2,4\mod 8$. This implies for $m\equiv 1,2,3,4,5\mod 8$ that~$\Sigma M$ can be obtained from the classical definition of the complex spinor bundle by forgetting the complex structure, and we will identify $\Sigma M$ with the classical spinor bundle in this case. In~the cases $m\equiv 2,3,4\mod 8$, there is even a quaternionic structure on $\Sigma_m$, again unique up to isomorphisms in the sense of quaternionic representations. Thus $\Sigma M$ is a quaternionic vector bundle. We write $\KK\ceq \HH$ for $m\equiv 2,3,4\mod 8$, $\KK\ceq \CC$ for $m\equiv 1,5\mod 8$, and $\KK\ceq \RR$ for~${m\equiv 0,6,7\mod 8}$.

The Dirac operator \smash{$\Dirac^g$} is an unbounded $\KK$-linear operator $\Gamma(\Sigma M)\to \Gamma(\Sigma M)$ with discrete spectrum. Thus all eigenspaces are (sub-)vector spaces over $\KK$.
The $\KO_m(\point)$-valued Atiyah--Singer index theorem says that
\[\alpha(M) =
 \begin{cases}
 \dim_\RR \ker\Dirac^g_+ -\dim_\RR \ker\Dirac^g_-\in \KO_m(\point)\cong \ZZ & \text{if}\ m\equiv0\mod 8,\\
 \dim_\CC \ker\Dirac^g\!\mod 2\in \KO_m(\point)\cong \ZZ/2 & \text{if}\ m\equiv 1\mod 8,\\
 \dim_\HH \ker\Dirac^g\!\mod 2\in \KO_m(\point)\cong \ZZ/2 & \text{if}\ m\equiv 2\mod 8,\\
 \dim_\HH \ker\Dirac^g_+ -\dim_\HH \ker\Dirac^g_-\in \KO_m(\point)\cong \ZZ & \text{if}\ m\equiv 4\mod 8,\\
 0\in \KO_m(\point)\cong \{0\} & \text{if}\ m\equiv 3,5,6,7\mod 8.
 \end{cases}
\]
The value of $\alpha(M)$ does not depend on the Riemannian metric and is thus a ``topological invariant''. More precisely, it only depends on the oriented homeomorphism type of $M$ for ${m\equiv 0,4 \mod 8}$, and it also depends on the differential structure and the spin structure if ${m\equiv 1,2\mod 8}$, but it is always independent of the Riemannian metric $g$.
In terms of the classical $\hat\maA$-genus, we have $\alpha(M)=\hat\maA(M)$ for $m\equiv 0\mod 8$ and $\alpha(M)=2\hat\maA(M)$ for $m\equiv 4\mod 8$.

In all cases, the index theorem provides a lower bound on $\dim_\KK\ker\Dirac^g$.
In fact, let us define
\[|\alpha(M)|\ceq
 \begin{cases}
 \bigl|\dim_\RR \ker\Dirac^g_+ -\dim_\RR \ker\Dirac^g_-\bigr| & \text{if}\ m\equiv0 \mod 8,\\
 1 & \text{if}\ m\equiv 1,2\mod 8 \ \text{and}\ \alpha(M)\neq 0,\\
 0 & \text{if}\ m\equiv 1,2\mod 8 \ \text{and}\ \alpha(M)= 0,\\
 \bigl|\dim_\HH \ker\Dirac^g_+ -\dim_\HH \ker\Dirac^g_-\bigr| & \text{if}\ m\equiv 4\mod 8,\\
 0 & \text{if}\ m\equiv 3,5,6,7,
 \end{cases}
\]
then the Atiyah--Singer index theorem implies $\dim_\KK\ker\Dirac^g\geq|\alpha(M)|$.
We define the space
\[
\Rmin(M) \ceq
\bigl\{g\in \maR(M) \mid \dim_\KK \ker \Dirac^{g}=|\alpha(M)| \bigr\}.
\]
of metrics for which the kernel of the Dirac operator attains this lower bound. In particular $\Rmin(M)=\Rinv(M)$ in the case $\alpha(M)=0$.

From now, one let us assume that~$M$ is connected. Then the set $\Rmin(M)$ is open and dense in $\maR(M)$, see \cite{ammann.dahl.humbert:09}. Metrics in $\Rmin(M)$ are called Dirac-minimal or \Dminimal{}, metrics in the complement $\maR(M) \setminus \Rmin(M)$ are called non-\Dminimal.

The main result of this article is the following theorem.
\begin{Theorem}\label{thm.main}
Let $M$ be a closed connected spin manifold of dimension $m=2$ or $m=4$. Then~$\Rmin(M)$ is connected.
\end{Theorem}

\begin{Remarks}\quad
\begin{enumerate}[(1)]\itemsep=0pt
\item
Note that a corresponding statement in many other dimensions does not hold. For example, for $\alpha(M)=0$ and $m\equiv 1,3,7,0\mod 8$, $m\geq 3$ the second author of this article \cite{dahl:2008} used the spectral flow for proving that $\Rinv(M)$ has more than one connected component.
This proof of disconnectedness cannot be extended to the case $\alpha(M)\neq 0$, however.
\item
Consider $m=2$ and let $\gamma$ be the genus of $M$. We will explain that non-\Dminimal{} Riemannian metrics exist if and only if
\begin{itemize}\itemsep=0pt
 \item $\alpha(M)=0$, $\gamma\geq 3$, or
 \item $\gamma\geq 5$.
\end{itemize}
Recall that Hitchin proved \cite[end of Section~2.1 and proof of Proposition~2.3]{hitchin:74} that $2\dim_\HH \ker\Dirac^g=\dim_\CC\ker\Dirac^g\leq \gamma+1$. Together with the $\KO(\point)$-valued index theorem, we conclude that all metrics are \Dminimal{} if $\gamma\leq 2$ or if $\alpha(M)\neq 0$ and $\gamma\leq 4$. In other words, non-\Dminimal{} metrics do not exist, when the above conditions are not satisfied.

Now let us assume $\gamma\geq 2$. Then $M$ carries a hyperelliptic complex structure, see \cite[Section~III.7]{farkas.kra:92} and we consider compatible Riemannian metrics~$g$.
In \cite{baer.schmutz:92}, the dimension of $\ker\Dirac^g$ was calculated in this setting by C.~B\"ar and P.~Schmutz-Schaller.
It follows from their calculation that for any integer $h$ with $0\leq 2 h \leq \gamma+1$ with $h\mod 2=\alpha(M)$ there is a~Riemannian metric~$g$ on $M$ with $\dim_\HH\ker\Dirac^g=h$. A Riemannian metric is non-\Dminimal{} if and only if $h\geq 2$. Thus non-\Dminimal{} metrics exist, if the above conditions are satisfied.

In some cases, a more geometric description can be obtained from the theory of minimal surfaces and the spinorial Weierstrass representation. We refer to work by Wei\ss{}, Witt and the first author \cite[Example~3.15 and Theorem~3.19]{ammann.weiss.witt.math.z:2016}, where the existence of metrics with nowhere vanishing harmonic spinors on $M$ is discussed. For $\gamma\geq 2$, such a metric necessarily has to be non-\Dminimal{}.
In the case $\alpha(M)= 0$ and $\gamma\geq 3$ and in the case $\alpha(M)\neq 0$, $\gamma\equiv 1\mod 4$ and $\gamma\geq 5$
embedded minimal surfaces in flat $3$-dimensional tori, diffeomorphic to~$M$, and with spin structure induced from this embedding exist. Their (spinorial) Weierstrass representation provides a space of harmonic spinors of quaternionic dimension $2$.
As the historical origin of this representation had attracted some interest, let us remark here that for the purpose of this remark the original literature by Weierstrass already gives a local description in terms of holomorphic functions which can be easily seen to glue together to a harmonic spinor.
\item If $M$ is a closed spin manifold of dimension $m$, and if $N$ is obtained from $M$ by a surgery of dimension $k\leq m-2$, then there is a Gromov--Lawson type surgery construction, that yields a map $\Rmin(M)\to \Rmin(N)$ \cite{ammann.dahl.humbert:09,ammann.humbert:2008,grosse.pederzani:2019}.
 In the spinorial context, the condition $k\leq m-2$ plays the same role as the condition $k\leq m-3$ does in the scalar curvature context.
 In the case $k\geq 1$, one can reverse the surgery, \ie $M$ is obtained from $N$ by a surgery of dimension $m-k-1\leq m-2$, and there is a similar map $\Rmin(N)\to \Rmin(M)$, and it is expected that this is a homotopy inverse to the map given above. Parts of this involved program were accomplished in \cite{grosse.pederzani:2019,pederzani:phd}.
\item So far, we have not yet commented on the regularity of the metrics in $\maR(M)$, nor on the topology of $\maR(M)$.
In our article, let $\ell\in \{1,2,3,\ldots\}\cup\{\infty\}$, and we consider $\maR(M)$ as the space of Riemannian metrics of regularity $C^\ell$, and we equip $\maR(M)$ with the $C^\ell$-topology.
We~write \smash{$\maR^{C^\ell}(M)$} for $\maR(M)$ when we want to emphasize the regularity.
For $\ell<\infty$, \smash{$\maR^{C^\ell}(M)$} is an open subset of the Banach space $C^\ell(M;T^*M\odot T^*M)$ and thus a Banach manifold. For $\ell=\infty$ it is a Fr\'echet manifold.
All our results hold for any $\ell$.
When proving Theorems~\ref{thm.main} and~\ref{thm.main.2}, we apply the submersion theorem in the Banach space setting and therefore restrict to $\ell<\infty$. We will discuss how to extend the results to $\ell=\infty$ in Remark~\ref{rem.regularity} and at the ends of the proofs of Theorems~\ref{thm.main.2} and~\ref{thm.main}.
\item In our discussion of scalar curvature in the introduction, we assume regularity $C^\ell$ with $\ell\geq 2$.
\item For a $C^1$-Riemannian metric~$g$, the Christoffel symbols in any chart are $C^0$, thus it makes sense to say that $g$ is \define{flat around} $p\in M$ if there is a $g$-parallel frame defined on a neighborhood of~$p$. This definition extends the classical definition of flatness for metrics of regularity~$C^2$.
A $C^1$-metric $g$ is called \define{conformally flat around} $p\in M$, if there is a $C^1$-function $u\colon M\to \RR_{>0}$ such that $ug$ is flat around $p$. We call $g$ (conformally) flat on~$M$, if it is (conformally) flat around any point of~$M$.
\end{enumerate}
\end{Remarks}

The result of Theorem~\ref{thm.main} can be refined to a result which essentially says that $\maR(M)\setminus\Rmin(M)$ is of codimension at least~$2$ in $\maR(M)$. In order to turn this into a precise statement, we have to introduce some notation.

For a subset $\maA$ of $\maR(M)$, we say that $\maA$ has
\define{codimension at least $k\in \NN\cup\{0,\infty\}$}
if for every $g\in \maA$ there is an open subset $\maU\subset\maR(M)$ of $g$ and a submanifold $\maS$ of $\maU\subset \maR(M)$ of codimension~$\geq k$ with $\maA\cap \maU\subset \maS$.
More precisely, $\maA$ has
codimension at least $k$, if every $g\in \maA$ has an open neighborhood $\maU$ and a submersion $F\colon\maU\to V$, $V$ open in $\RR^k$, with $F(\mathcal{A}\cap\mathcal{U})=\{0\}$.

In the remaining part of the introduction, we assume $m\in\{2,4\}$.
We define the space $\Rcf(M)$ as follows. If~$M$ is diffeomorphic to $T^4$ we define
\[
\Rcf(M)\ceq \bigl\{u\cdot g\mid u\in C^\ell(M;\RR_{>0}) \text{ and $g$ is a flat metric}\bigr\}
\]
and $\Rcf(M)=\varnothing$ if~$M$ is not diffeomorphic to $T^4$. In particular, $\Rcf(M)=\varnothing$ in the case $m=2$.
The set $\Rcf(M)$ is a closed subset of $\maR(M)$ of infinite codimension, see Proposition~\ref{prop.CF.codim}.

Note that $\ker\Dirac^g$ is a quaternionic vector space. In the following, we write $\omega\in \Cl(TM)$ for the \define{real volume element}, \ie if $(e_1,\ldots,e_m)$ is a positively oriented orthonormal basis of some $(T_xM,g_x)$, then $\omega|_{x}=e_1\cliff\cdots\cliff e_m$.
The action of~$\omega$ on spinors anticommutes with~$\Dirac$ as~$m$ is even.

For $m=4$, we have $\omega\cliff \omega=1$ and thus the action of~$\omega$ on spinors induces a decomposition $\Sigma M=\Sigma_+M\oplus \Sigma_- M$ into quaternionic subspaces for the eigenvalues $\pm 1$ of the idempotent $\omega$. Be aware that the standard conventions imply that $\omega=-1$ on $\Sigma_+M$ and $\omega=+1$ on $\Sigma_-M$.
This yields a decomposition \smash{$\ker\Dirac^g=\ker\Dirac^g_+\oplus \ker\Dirac^g_-$} and we have
\[
\hat \maA(M)/2 = \dim_\HH \ker\Dirac^g_+ - \dim_\HH \ker\Dirac^g_-.
\]
Thus, using \smash{$|\alpha(M)|\ceq \bigl| \hat \maA(M)/2\bigr|$} we see that for $g\in \maR(M)$ and
\[
d\ceq\min\bigl\{\dim_\HH \ker\Dirac^g_+,\dim_\HH \ker\Dirac^g_-\bigr\}\in \NN_0,
\]
we have \smash{$\dim_\HH \ker\Dirac^g=|\alpha(M)| +2d$}.

For $m=2$, we have $\omega\cliff \omega=-1$, and thus the map $\phi\mapsto \omega\cliff \phi  i $ is an idempotent as well. Its associated decomposition $\ker\Dirac^g=\ker\Dirac^g_+\oplus \ker\Dirac^g_-$ is however no longer invariant under multiplication with $j$. In fact, one easily checks that multiplication with $j$ defines antilinear isomorphisms $\ker\Dirac^g_+\to\ker\Dirac^g_-$ and $\ker\Dirac^g_-\to\ker\Dirac^g_+$. Then ${\dim_\HH\ker\Dirac^g=\dim_\CC\ker\Dirac^g_+=\dim_\CC\ker\Dirac^g_-}$ and it turns out that $\dim_\HH\ker\Dirac^g\mod 2$ can be identified with the index $\alpha(M)$ of $M$ in $\KO_2(\point)\cong\Z/2\Z$. We write $|\alpha(M)|=1$ if this index is nontrivial and $|\alpha(M)|=0$ if it is trivial. Then there is a number $d\in \NN_0$, such that $\dim_\HH\ker\Dirac^g=|\alpha(M)|+2d$.

In both cases $m=2$ and $m=4$, we define
\[
\maR_d(M)\ceq \bigl\{g\in \maR(M)\setminus\Rcf(M)\mid \dim_\HH \ker \Dirac^{g}=|\alpha(M)|+2d\bigr\}.
\]
In particular, $\Rmin(M)=\maR_0(M)$ and we have the disjoint union
\[
\maR(M)=\Rcf\amalg \coprod_{d\in \NN_0}\maR_d(M).
\]
Note that \smash{$\maR_{\leq d}(M) \ceq \bigcup_{0\leq l\leq d}\maR_l(M)$} is an open subset of $\maR(M)$.

\begin{Theorem} \label{thm.main.2}
Let $M$ be a closed connected spin manifold of dimension $m=2$ or $m=4$. Then for $d\in \NN$, $d>0$, $\maR_d(M)$ has codimension at least~$2$ in $\maR_{\leq d}(M)$. Furthermore, $\Rcf(M)$ has infinite codimension in $\maR(M)$.
\end{Theorem}

Note that we do not expect the codimension estimate for $\maR_d(M)$ to be sharp except in the case $m=2$, $d=1$ and in the case that we have $m=4$, $d=1$ and that~$M$ is ``almost'' a product of surfaces.
Improved estimates in the other cases shall be the subject of a future article.
However, the codimension $2$-result is what we need for the connectedness of $\Rmin(M)$ which is the topic of this article.

For the following investigations, we expect that the following statement will be helpful.
This theorem is also an important intermediate step for the proof of Theorem~\ref{thm.main.2}, but it allows finer information about the behavior of harmonic spinors deforming into eigenspinors with small eigenvalues.
It does not only provide information about the eigenvalues, but also about the way a kernel may break up in eigenspaces for different small eigenvalues.

In the following, we denote the $\HH$-valued $L^2$-scalar product on $L^2(\Sigma^gM)$ by \smash{$\llangle\biargu\rrangle_\HH$}.

\begin{Theorem} \label{thm.main.3}
 Let $(M,g_0)$ be a closed connected Riemannian spin manifold of dimension $m=2$ or $m=4$. In the case $m=4$, we additionally assume $g_0\notin \Rcf(M)$.
 Let $B_\epsilon(0)$ be the open \ball{} of radius $\epsilon$ around $0$ in $\RR^2$.

Suppose that we have non-zero $\Phi\in \ker \Dirac_+^{g_0}$ and $\Psi\in \ker\Dirac_-^{g_0}$. In the case $m=2$, we additionally assume that $\llangle\Phi,\Psi\rrangle_\HH=0$. Then there is some $\epsilon>0$ and a smooth family $g_{\bullette}\colon B_\epsilon(0)\to \maR(M)$ with $g_{s,t}\in \maR(M)$, $g_{0,0}=g_0$ with the property:
For any smooth extensions $\Phi_{s,t}\in \Gamma(\Sigma_+^{g_{s,t}}M)$, $\Psi_{s,t}\in \Gamma(\Sigma_-^{g_{s,t}}M)$, $(s,t)\in B_\epsilon(0)$, with $\Phi_{0,0}=\Phi$ and $\Psi_{0,0}=\Psi$, we have for any choice of $(s,t)\in \RR^2\setminus \{(0,0)\}$:
\[
\frac{\upd}{\upd \tau}\bigg|_{\tau=0}
\llangle\Dirac^{g_{\tau s,\tau t}}\Phi_{\tau s,\tau t} \Psi_{\tau s,\tau t}\rrangle_\HH
\neq 0 .
\]
\end{Theorem}

Note that the condition $\llangle\Phi,\Psi\rrangle_\HH=0$ is automatically satisfied for $m=4$, but not for $m=2$, see Section~\ref{subsec.quaternionic-structures}.

Theorem~\ref{thm.main.3} will be important for two reasons. At first, it provides the essential step in the proof of Theorem~\ref{thm.main.2}.
Secondly, it will allow further applications in work in progress. In~particular, Theorem~\ref{thm.main.3} implies the following two statements: If $M$ is a closed connected spin surface
with $\alpha(M)\neq 0$, then the space of non-minimal metrics is of codimension at least $4$ in the space of all Riemannian metrics. If $M$ is a closed connected $4$-dimensional spin manifold with $\bigl|\hat A(M)\bigr|=2k$, then the space of non-minimal metrics is of codimension at least $2(k+1)$ in the space of all Riemannian metrics. We will elaborate on this in a forthcoming article.

The structure of the article is as follows. In Section~\ref{sec.prelim}, we collect some preliminaries and introduce some notation and some conventions used in the article. In particular, we recall work by Bourguignon and Gauduchon \cite{bourguignon.gauduchon:92} that gives the variational formula for the Dirac operator under variations of the metric. The essential term in the associated variational formula is the energy-momentum tensor which will lead to further investigations in Section~\ref{sec.def.harm.spinors}. Here we analyze the energy-momentum tensor at a non-\Dminimal{} metric, and we will find a formula, namely equation~\eqref{eq.GZ-prod} that establishes a relationship between the energy-momentum tensor associated to a pair of harmonic spinors $(\Phi,\Psi)$ and the gradient of $|\Phi|^2/|\Psi|^2$. This formula will be used in Section~\ref{sec.proof.thm.main.2} for proving Theorem~\ref{thm.main.3}, and then we deduce Theorem~\ref{thm.main.2}. The proof of Theorem~\ref{thm.main} is then given in Section~\ref{sec.proof.thm.main}.

\section{Notation and preliminaries}\label{sec.prelim}

\subsection{Scalar products on quaternionic vector spaces}

See also \cite[Chapter~2]{harvey:90}.

We follow the convention that vectors in quaternionic vector spaces are multiplied by quaternionic scalars from right. However we allow real scalars to be multiplied both from the left and from the right.

Let $V$ be a quaternionic vector space. An $\RR$-bilinear map $g\colon V\times V \to \HH$ is \define{$\HH$-sesquilinear} if we have
\begin{equation*}
g(v\lambda,w)= \bar{\lambda} \cdot g(v,w),\qquad g(v,w\lambda)= g(v,w)\lambda
\end{equation*}
for every $v,w\in V$ and every $\lambda\in \HH$. We then also say that $g$ is an \define{$\HH$-sesquilinear form}.
Such a~form is called \define{$\HH$-hermitian} if for every $v,w\in V$
\begin{equation}\label{eq.HH.hermitian}
g(w,v)= \overline{g(v,w)},
\end{equation}
and \define{positive definite} if for every $v\in V\setminus \{0\}$ we have $g(v,v)>0$.

\begin{Lemma}
A positive definite $\HH$-sesquilinear form is $\HH$-hermitian.
\end{Lemma}

\begin{proof}
Let $I\in \spann\{i,j,k\}$ be one of the complex structures in $\HH$, \ie $I^2=-1$.
Let $\pi_I\colon \HH\to \CC$ be the orthogonal projection with $\pi_I(1)=1$, $\pi_I(I)=i$.
Then $\pi_I\circ g$ is a positive definite sesquilinear form on $(V,I)$ in the complex sense.
From classical linear algebra over $\CC$ we know that $\pi_I\circ g$ is hermitian, \ie \eqref{eq.HH.hermitian} is satisfied for $\pi_I\circ g$ instead of $g$.
Applying this for $I\in \{i,j,k\}$ implies \eqref{eq.HH.hermitian}.
 \end{proof}

An \define{$\HH$-scalar product} or a \define{quaternionic scalar product} is defined as a positive definite $\HH$-sesqui\-linear form.

\begin{Example}
We write \smash{$X^*=\overline{X^{\mathsf T}}= \ol{X}^{\mathsf T}$} for $X\in \HH^n$. On $\HH^n$, we have the standard scalar product
\[
\<X,Y\>_\HH= X^* Y.
\]
\end{Example}
If $\<\biargu\>_\HH$ is a quaternionic scalar product on $V$, then $\<\biargu\>\ceq \Re \<\biargu\>_\HH$ is a real scalar product and it satisfies
\begin{equation}\label{eq.HH.linearity}
\<\phi,\psi\>_\HH = \<\phi,\psi\>- \<\phi,\psi\cdot i \> \cdot i - \<\phi,\psi\cdot j \> \cdot j - \<\phi,\psi\cdot k \> \cdot k
\end{equation}
for any $\phi,\psi\in V$. Further, multiplications by $ i $, $ j $ and $ k $ are isometries with respect to
$\<\biargu\>$.

Conversely, if $\<\biargu\>$ is a real scalar product such that $ i $, $ j $ and $ k $ are isometries, then we may use \eqref{eq.HH.linearity} to define a quaternionic scalar product with $\<\biargu\>= \Re \<\biargu\>_\HH$.
In particular, $\<\biargu\>_\HH$ is $\HH$-sesquilinear.

Note that in both settings we then also have
\begin{equation*} 
\<\phi,\psi\>_\HH = \<\phi,\psi\>+\<\phi\cdot i ,\psi\> \cdot i + \<\phi\cdot j ,\psi\> \cdot j + \<\phi\cdot k ,\psi\> \cdot k .
\end{equation*}

Let $V$ be a quaternionic line.
For $\psi\in V$ with $|\psi|=1$, one checks that $(\psi,\psi\cdot i ,\psi \cdot j , \psi\cdot k )$ is an orthonormal real basis of the quaternionic line $V$. As a consequence, we get
\begin{align}
|\psi|^2\phi
=
\psi \cdot \<\psi,\phi\>
+ \psi \cdot i \<\psi\cdot i ,\phi \>
+ \psi \cdot j \<\psi\cdot j ,\phi \>
+ \psi \cdot k \<\psi\cdot k ,\phi\>
=
\psi \cdot \<\psi,\phi\>_\HH\label{formula.HH.basis}
\end{align}
for any $\phi,\psi\in V$.

\subsection{Spinors and quaternionic structure}\label{subsec.quaternionic-structures}

We want to recall some facts about spin structures and spinors on oriented manifolds $M$ of dimension $m$.

In the literature several equivalent definitions of spin structures may be found.
In differential geometry one often defines a spin structure by considering the principal $\SO(m)$-bundle $P_\SO(M)$ of positively oriented orthonormal frames, and then one defines a spin structure as an equivariant double cover $P_\Spin(M) \to P_\SO(M)$ \cite{friedrich-book,lawson.michelsohn:89}.
This type of spin structure, called a \define{metric spin structure}, has the drawback that it depends on the choice of a Riemannian metric~$g$.
An~alternative way is to consider the non-trivial double cover \smash{$\witi\GL_+(m,\RR)\to \GL_+(m,\RR)$} instead and to define a spin structure as an equivariant double cover \smash{$P_{\witi\GL_+}(M) \to P_{\GL_+}(M)$}, where $P_{\GL_+}(M)$ is the principal \smash{$\GL_+(m,\RR)$}-bundle of positively oriented frames. This type of definition allows us to define a spin structure on an oriented manifold, in a way independent of a Riemannian metric, and will be called a \define{topological spin structure} or later simply a \define{spin structure}.
It is classically known since many decades that there is a 1-to-1 relation between isomorphism classes of metric spin structures and isomorphism classes of topological spin structures, but as this occasionally led to misunderstandings, it is helpful that details of this construction were worked out, \eg in~\cite{swift-parttwo} and~\cite[Section~2.2]{ammann.humbert:2008}.
Let us also briefly mention that in topology spin structures are often defined as a lift of the classifying map $M\to B\SO$ for $TM$ to a map $M\to B\Spin$ whose homotopy classes again provide an equivalent definition; this point of view will not be used in this article.

Thus we will assume that $M$ carries a (topological) spin structure in this sense, and that it defines -- for each Riemannian metric~$g$ -- a spin structure \smash{$P_\Spin^g(M) \to P_\SO^g(M)$} in the usual Riemannian sense. Using a representation
$\rho\colon\Spin(m)\to \GL(V)$ on a vector space $V$,
we define an associated bundle \smash{$V^gM\ceq P_\Spin^g(M) \times_\rho V$}.

We will use this construction for spinor bundles, where $V= \Sigma_m$ is an irreducible linear Clifford module.
An \define{irreducible linear Clifford module} is a representation $\rho\colon\Spin(m)\to \GL(\Sigma_m)$ that extends to an irreducible representation of the Clifford algebra $\Cl(\RR^m)$.
This definition comes in different flavors, in particular one may choose different structures on $\Sigma_m$. For example, one may consider $\Sigma_m$ as a real vector space, and then one may show that $\Sigma_m$ has a scalar product for which Clifford-multiplication is skew and thus $\rho$ acts isometrically.
The scalar product on $\Sigma_m$ is then uniquely determined up to scaling, thus we obtain an equivalent definition if we view~$\Sigma_m$ as a Euclidean vector space and if we consider representation $\rho\colon\Spin(m)\to \O(\Sigma_m)$ instead. Similarly, connectedness of $\Spin(m)$ allows to reduce to representations $\rho\colon\Spin(m)\to \SO(\Sigma_m)$.

For each representation $\rho\colon\Cl(\RR^m)\to \End(\Sigma_m)$, the associated bundle construction provides an associated \define{spinor bundle} $\Sigma^{g,\rho} M\ceq P^g_\Spin(M)\times_\rho\Sigma_m$. The $\Spin(m)$-equivariance of the map~$\rho$ leads to an action of $\Cl(TM,g)$ on $\Sigma^{g,\rho} M$, see again \cite{friedrich-book,lawson.michelsohn:89} for details. For $p\in M$, the action of $X\in \Cl(T_pM,g)$ on \smash{$\phi\in \Sigma_p^{g,\rho} M$} is denoted by $X\cliff \phi$, and we also use $X\cliff Y$ to denote the multiplication of $X,Y\in\Cl(TM,g)$.

For a given Riemannian metric~$g$, we obtain the real spinor bundle $\Sigma^g_\RR M\to M$.
Note that the isomorphism type of the real representation $\rho\colon\Spin(M)\to \GL(\Sigma_m)$, of real orthogonal representation $\rho\colon\Spin(M)\to \O(\Sigma_m)$ is not unique for $m\equiv 3\mod 4$, but it is unique for our main cases of interest, namely $m\in \{2,4\}$. We refer to \cite[Chapter~I, Section~5]{lawson.michelsohn:89} for a discussion of such representations.

Similarly, we may require that $\Sigma_m$ is a vector space over the complex numbers~$\CC$ or over the quaternionic numbers~$\HH$,
and we obtain \smash{$\Sigma^g_\CC M$} and \smash{$\Sigma^g_\HH M$}
in the same way by replacing real, bilinear scalar products by complex or quaternionic sesquilinear scalar products.
We denote the resulting real, complex or quaternionic fiberwise scalar product as $\<\biargu\rangleR$, $\<\biargu\rangleC$, or $\<\biargu\rangleH$.
By integrating with respect to the Riemannian volume element, we obtain for $\KK\in \{\RR,\CC,\HH\}$
the $\KK$-sesquilinear $L^2$-scalar product on the space of $L^2$-sections $L^2\bigl(\Sigma^g_\KK M\bigr)$, denoted by $\llangle\biargu\rrangle_\KK$. In~the notation of formulas, we suppress the Riemannian volume element in integrals for the sake of readability.

Now we specialize to the cases $m\equiv 2,4\mod 8$ in which there is up to isomorphism one irreducible real (resp.\ complex, resp.\ quaternionic) representation $\Sigma^\RR_m$ \big(resp.\ $\Sigma^\CC_m$, resp.\ $\Sigma^\HH_m$\big) of~$\Cl(\RR^m)$, see \cite[Chapter~I, Section~5]{lawson.michelsohn:89} and \cite[Section~1.7]{friedrich-book}. As these spaces have the same real dimension, there is a complex isomorphism $\Sigma^\HH_m\to \Sigma^\CC_m$ and a real isomorphism $\Sigma^\CC_m\to \Sigma^\RR_m$. When passing to the associated bundles, we obtain isomorphisms $\Sigma^g_\HH M\to \Sigma^g_\CC M\to \Sigma^g_\RR M$,
``forgetting'' multiplication with $ j $ and then with $ i $, which we use to identify these bundles.

The \define{real volume element} in the Clifford algebra is defined as $\omega\ceq e_1\cliff e_2\cliff \cdots\cliff e_m\in \Cl(T_pM)$, for a positively oriented orthonormal basis $(e_a)$. For $m$ even, the \define{complex volume element} is defined as $\omega_\CC\ceq i ^{m/2} \omega \in \Cl(T_pM)\otimes_\RR\CC$, for a positively oriented orthonormal basis $(e_a)$. Because of $\omega_\CC^2=1$, it defines a grading $\Sigma^g_\CC M = \Sigma^g_+ M\oplus \Sigma^g_-M$ in the sense of complex bundles with scalar product.

For $m=4$, quaternionic multiplication commutes with $\omega_\CC$, and thus this splitting also holds in the sense of quaternionic bundles with $\HH$-sesquilinear scalar product. In particular, this implies that $\llangle\Phi,\Psi\rrangle_\HH=0$ for $\Phi\in \Gamma\bigl(\Sigma^g_+ M\bigr)$ and $\Psi\in \Gamma(\Sigma^g_- M)$.

For $m=2$, quaternionic multiplication by $ j $ anti-commutes with $\omega_\CC$; thus $\Sigma^g_+ M$ and $\Sigma^g_-M$ are no longer quaternionic bundles, and multiplication by $ j $ defines complex anti-linear isomorphisms~${\Sigma^g_+ M\to \Sigma^g_-M}$ and $\Sigma^g_- M\to \Sigma^g_+M$. In this dimension $\Phi\in \Gamma\bigl(\Sigma^g_+ M\bigr)$ and $\Psi\in \Gamma(\Sigma^g_- M)$ satisfy $\llangle\Phi,\Psi\rrangle_\CC=0$, but, \eg for $\Psi\ceq \Phi\cdot j $ we have $\llangle\Phi,\Psi\rrangle_\HH= \|\Phi\|^2\cdot j $.

\subsection{The Dirac operator}

The Riemannian metric $g$ on $M$ gives the Levi-Civita connection $\nabla^g$ on $TM$, and an associated connection-1-form on \smash{$P_\SO^g(M)$}. This connection-1-form lifts to a connection-1-form on any spin structure \smash{$P_\Spin^g(M)$} and finally induces a connection, denoted by $\nabla^g$ as well, on any associated spinor bundle $\Sigma^g M$, see \cite[Chaper~II, Section~4]{lawson.michelsohn:89} and \cite[Section~3.1]{friedrich-book}.

The \define{Dirac operator} $\Dirac^g\colon \Gamma(\Sigma^gM)\to \Gamma(\Sigma^gM)$ is defined by
\[
\Dirac^g \phi \ceq \sum_{a=1}^n e_a \cliff \nabla^g_{e_a} \phi,
\]
where $(e_a)$ is a locally defined orthonormal frame.
The Dirac operator is a self-adjoint elliptic densely defined unbounded operator $\Dirac^g \colon L^2(\Sigma^g M) \to L^2(\Sigma^g M)$.
As such, since~$M$ is closed, the Dirac operator is a Fredholm operator with unbounded real spectrum consisting of eigenvalues with finite-dimensional eigenspaces, see \cite[Chapter~III, Section~5]{lawson.michelsohn:89} and \cite[Section~4]{friedrich-book}. In particular, the kernel $\ker\Dirac^g$ is a finite-dimensional vector space. A \define{harmonic spinor} is an element of $\ker\Dirac^g$.

When we use the $C^1$-topology on $\maR(M)$, any metric has a neighborhood $\maU$ such that the eigenvalues
of $\Dirac^g$ for $g\in \maU$ can be numbered by functions that are continuous in $g$, see, for example, \cite{nowaczyk:continuity}.

\subsection[Identifying spinor bundles and Dirac operators over spaces of Riemannian metrics]{Identifying spinor bundles and Dirac operators\\ over spaces of Riemannian metrics} 

On a spin manifold $M$ with a fixed (topological) spin structure,
not only the connection and the Dirac operator but even the construction of the spinor bundle
depends on the choice of a Riemannian metric on $M$ in a subtle way. The dependence of the choice of metric was investigated by Bourguignon and Gauduchon in \cite{bourguignon.gauduchon:92}.
Given two metrics $g$, $h$, there is an endomorphism \smash{$b^g_h$} of $TM$
such that \smash{$h\bigl(b^g_h \argu, b^g_h \argu\bigr) = g(\argu,\argu)$}, and \smash{$b^g_h$} is symmetric
and positive definite (with respect to both metrics).
The endomorphism $b^g_h$ is uniquely determined by these properties and
gives a~map \smash{$P_\SO^g(M) \to P_\SO^h(M)$}.
Further, this map lifts to a map \smash{$\beta^g_h \colon P_\Spin^g(M) \to P_\Spin^h(M)$} between the $\Spin(m)$-principal bundles defined by the same topological spin structure.
This, in turn induces linear maps
\[
\beta^g_h\colon\ \Sigma^g M\to \Sigma^h M, \qquad
\beta^g_h\colon\ \Gamma(\Sigma^g M) \to \Gamma\bigl(\Sigma^h M\bigr)
\]
of the associated spinor bundles and their spaces of smooth sections.
The map~\smash{$\beta^g_h$} is by construction a fiberwise isometry of the Hermitian metrics. We recalibrate this map by considering the volume elements $\dV^{g}$ and \smash{$\dV^{h} = \bigl(f^g_h\bigr)^2\dV^{g}$} of~$g$ and~$h$, where \smash{$f^g_h$} is a positive function.
Following~\cite{maier:97}, we define the rescaled map \smash{$\hat\beta^g_h= \bigl(f^g_h\bigr)^{-1} \beta^g_h$},
and we get another vector bundle isomorphism \smash{$\hat\beta^g_h\colon \Sigma^g M\to \Sigma^h M$}, which is not a fiberwise isometry, and an isometry
\[
\hat\beta^g_h\colon\ L^2(\Sigma^g M) \to L^2\bigl(\Sigma^h M\bigr)
\]
of the corresponding Hilbert spaces of $L^2$-sections. We have \smash{$\beta^h_g\circ \beta^g_h=\id$} and \smash{$\hat\beta^h_g\circ \hat\beta^g_h=\id$}.
However, one has to be aware of the effect that in general for three Riemannian metrics~$g$,~$h$,~$k$ we have \smash{$\beta^h_k\circ \beta^g_h\neq\beta^g_k$} and \smash{$\hat\beta^h_k\circ \hat\beta^g_h\neq\hat\beta^g_k$}, so the identification given by \smash{$\hat\beta^g_h$} is not transitive.

Having fixed a Riemannian metric $g_0$, the map \smash{$g\mapsto \hat\beta_{g_0}^g$} defines trivializations
\[
\Sigma^{\univ}_p M\ceq \coprod_{g\in \Sym_+T^*_pM }\Sigma^g_p M\to \Sigma^{g_0}_p M\times \Sym_+T^*_pM
\]
for any $p\in M$ and
\[
L\ceq \coprod_{g\in \maR(M)}L^2(\Sigma^g M)\to L^2(\Sigma^{g_0}M)\times \maR(M),
\]
and the latter map is a fiberwise isometry of Hilbert spaces from the bundle $L\to \maR(M)$ to the trivial $L^2(\Sigma^{g_0}M)$-bundle over $\maR(M)$.
Because of the lacking transitivity for $\hat\beta^g_h$, this identification depends on the choice of~$g_0$. This way to trivialize the Hilbert bundle $L\to \maR(M)$ is worked out by Maier \cite{maier:97}, relying on \cite{bourguignon.gauduchon:92}.

For understanding the results of this article, it is sufficient to work in Maier's trivialization.
However, let us briefly sketch in the next subsection a more conceptual picture,
which may turn out to be important when considering, for instance,
higher order terms in a perturbative expansion.

\subsection{The universal spinor bundle}

In this subsection, we describe a natural connection $\hat\nabla$ on the Hilbert bundle $L\to \maR(M)$. As~the concepts of this sections are not required for the logical structure of the article, we will only sketch this connection. However this connection will motivate our notation used later on.

In the following, let $\Sigma^g M\to M$ be the real, complex or quaternionic spinor bundle over a~Riemannian spin manifold $(M,g)$.
Note that for $p\in M$, the fiber $\Sigma^g_p M$, viewed as a Clifford module with compatible scalar product only depends on $g_p\in \Sym_+T^*_pM$ and thus one may define the \define{universal spinor bundle}
\[\Sigma^{\univ} M\ceq\coprod_{p\in M}\coprod_{g\in \Sym_+T^*_pM }\Sigma^g_p M\to \Sym_+ T^*M ,\]
see \cite{ammann.weiss.witt.math.ann:2016,ammann.weiss.witt.math.z:2016,mueller.nowaczyk:2017} for details.
This bundle carries a natural partial connection, where ``partial'' means that only differentiation along the fibers of the bundle $\Sym_+ T^*M\to M$ is defined. Now, the Hilbert bundle $L\to \maR(M)$ inherits a connection $\nabla$ from the partial connection on~$\Sigma^{\univ} M$. This connection is tightly related to the maps \smash{$\beta^g_h$} discussed above from \cite{bourguignon.gauduchon:92}.
On the one hand, if we have $\phi\in \Gamma(L)$ \big(\ie $\phi_h\in L^2(\Sigma^h M)$, depending differentiably on $h\in \maR(M)$\big) and ${k\in T_g\maR(M)\cong \Gamma\big(\Sym T^*M\big)}$, then we have the formula
\[
\nabla_k\phi=\frac{{\rm d}}{{\rm d}t}\bigg|_{t=0}\beta^{g+tk}_g\phi_{g+tk} .
\]
On the other hand, $\beta^g_h$ is the parallel transport along $t\mapsto (1-t)g+th$ with respect to~$\nabla$.
Note that this parallel transport does \emph{not} preserve the $L^2$-scalar product on~$L$,
thus $\nabla$ is not compatible with the $L^2$-metric.

Recall from the preceding subsection, that Maier recalibrated the map~\smash{$\beta^g_h$} to~\smash{$\hat\beta^g_h$} in order to preserve the fiberwise $L^2$-scalar product $\llangle\biargu\rrangle$ on~$L\to \maR(M)$.
In the same spirit, we define a~\define{recalibrated connection} $\hat\nabla$ on $L$ by the formula
\[\hat\nabla_k\phi=\frac{{\rm d}}{{\rm d}t}\bigg|_{t=0}\hat\beta^{g+tk}_g\phi_{g+tk} .\]
This connection also arises from a recalibrated connection on
 $\Sigma^{\univ} M\to \Sym_+ T^*M$.
 As \smash{$\hat\beta^g_h$} preserves the $L^2$-scalar product, the connection $\hat\nabla$ is compatible with the fiberwise $L^2$-metric on $L\to\maR(M)$, \ie
\[\partial_k \llangle\phi,\psi\rrangle= \bigl\langle\!\!\bigl\langle\hat\nabla_k \phi,\psi\bigr\rangle\!\!\bigr\rangle+ \bigl\langle\!\!\bigl\langle\phi,\hat\nabla_k \psi\bigr\rangle\!\!\bigr\rangle\]
for $\phi,\psi\in L^2(\Sigma^g M)$ and $k\in \Gamma\bigl(\Sym T^*M\bigr)\cong T_g\maR(M)$.

This construction may be done for real, complex and quaternionic spinor bundles
and we get forgetful isomorphisms
\[\Sigma^{\univ}_\HH M\to\Sigma^{\univ}_\CC M\to \Sigma^{\univ}_\RR M\]
by forgetting multiplication with $j$ and $i$, respectively. These forgetful isomorphisms preserve the connections $\nabla$ and $\hat\nabla$ and the fiberwise norm. In dimension~$2$ and~$4$, these are isomorphisms of complex vector bundles and real vector bundles, respectively.

At the end of this excursion on the universal spinor bundle,
let us mention that for a given background metric $g_0\in \maR(M)$ we have four connections on the bundle $L \to \maR(M)$.
The first is defined by the Bourguignon--Gauduchon trivialization~$\beta_{g_0}^g$, the second by Maier's trivialization~$\hat\beta_{g_0}^g$, the third and the fourth are the connections~$\nabla$ and~$\hat\nabla$.
We will prove a transversality statement for a section of a related bundle over $\maR(M)$ at a metric $g_0$ where the section vanishes.
This transversality statement is unaffected by the choice of connection.

\subsection{Variation of the Dirac operator with respect to the metric} \label{subsec:variation.of.D}

From now, we specialize to dimension $m=\dim M\in \{2,4\}$ for which the spinor bundles are quaternionic line bundles.

For the bundle $L\to\maR(M)$ defined above, let $\End(L)\to \maR(M)$ be the bundle of densely defined fiberwise $\mathbb{H}$-linear operators.
The Dirac operator $\Dirac^{\bullette}\colon g \mapsto \Dirac^g$ is a section of this bundle, and due to its symmetries it is even a section of the subbundle \smash{$\End_{\omega\text{-ac}}^{\text{SA}}(L)\to\maR(M)$} of self-adjoint operators anti-commuting with the real volume element $\omega$.

It is a classical fact, see \cite[Proposition~2.4]{maier:97} for a reference, that $\Dirac^{ \bullette}$ is differentiable in the sense that we may derive $\Dirac^{ \bullette}$ with respect to the induced connection $\hat\nabla$ in any direction $k\in T_{g_0}\maR(M)\cong\Gamma\bigl(\Sym T^*M\bigr)$, and the derivative is given by
\begin{equation*} 
 \frac{{\rm d}}{{\rm d}t}\bigg|_{t=0}\bigl(\hat\beta^{g_0+tk}_{g_0}\circ \Dirac^{g_0+tk}\circ\hat\beta_{g_0+tk}^{g_0}\bigr)=
-\frac{1}{2} \sum_{a,b=1}^m k(e_a,e_b) e_a \cliff \nabla^{g_0}_{e_b}
-\frac{1}{4} (\div^{g_0} k) \cliff ,
\end{equation*}
where $(e_1,\dots,e_m)$ is a locally defined orthonormal frame, for $k\in \Gamma\big(\Sym T^*M\big)$.
In view of~the interpretation of \smash{$\hat\beta^g_h$} as the parallel transport for the connection \smash{$\hat\nabla$}, we write the left-hand~side of this equation as
\smash{$\bigl(\frac{\hat\nabla}{{\rm d}g}\big|_{g=g_0} \Dirac^{g}\bigr) (k)$},
so we obtain
\begin{equation} \label{first-var-D}
\biggl(\frac{\hat\nabla}{{\rm d}g}\bigg|_{g=g_0} \Dirac^{g}\biggr) (k)
=
-\frac{1}{2} \sum_{a,b=1}^m k(e_a,e_b) e_a \cliff \nabla^{g_0}_{e_b}
-\frac{1}{4} \div^{g_0} k \cliff .
\end{equation}

The locally defined orthonormal frame $(e_1,\dots,e_m)$ allows taking traces by contraction, \eg for $A\in \End(TM)$ and $\alpha,\beta\in \Omega^1(M)$ we have
\[
\tr(A)=\sum_{a=1}^mg(A(e_a),e_a),\qquad \tr^g(\alpha\otimes \beta)=\<\alpha,\beta\rangleR=\sum_{a=1}^m\alpha(e_a)\beta(e_a) .
\]
In~particular, the right-hand side does not depend on the choice of the frame. However, the right-hand side of these equations is only defined on the domain of the corresponding frame. In~the following, we abuse notation by interpreting the formulas on the right side of \eqref{first-var-D} as formulas on~$M$.
The above expression then means that for any $p\in M$ and any frame $(e_1,\dots,e_m)$ defined on an open neighborhood $U_p$ of $p$, the formula holds on~$U_p$. Obviously, in order to check this, it is sufficient to prove that for any $p\in M$ this holds for some frame defined on some $U_p\ni p$.

\subsection{Codimension of conformally flat metrics}

\begin{Proposition} \label{prop.CF.codim}
Let $M$ be a manifold of dimension at least $4$. Then the set of locally conformally flat metrics has infinite codimension in $\maR(M)$.
\end{Proposition}

\begin{proof}
 In dimensions at least 4, a Riemannian metric is locally conformally flat if and only if its Weyl curvature tensor vanishes identically.
 This fact, which is usually formulated for smooth metrics, easily extends to metrics of regularity $C^\ell$ for $\ell\geq 2$, an even holds in some weak sense for $\ell=1$.

Let $g^0$ be locally conformally flat. Let $ p \in M$. In a neighborhood~$U$ of~$p$ there are coordinates~$x^a$ centered at $p$ in which $g^0_{ab} = u \delta_{ab}$ for some positive $C^\ell$-function $u$. Let $w_{acdb}$ be a~constant tensor with all the pointwise symmetries of a Weyl tensor, and set
\[
\gamma^t_{ab} \ceq \delta_{ab} - \frac{1}{3} t \cdot \chi \sum_{c,d=1}^n w_{acdb} x^c x^d,
\]
where $\chi$ is a cut-off function with support in $U$ and $\chi=1$ near $p$. Set $g^t_{ab} = u \gamma^t_{ab}$. Then the $(4,0)$-Weyl curvature tensor of $g^t$ is $W^{g^t} = u W^{\gamma^t}$ and at the center of coordinates we have \smash{$W^{\gamma^t}(p)_{acdb} = tw_{acdb}$}. Let $\maW_n$ be the dimension of the space of pointwise Weyl tensors in dimension~$n$. Then the above construction provides a $\maW_n$-dimensional space of perturbations~$g^t$ which deform $g^0$ away from being conformally flat, so the space of locally conformally flat metrics has codimension at least $\maW_n$. Repeating the argument with any finite number points $p$ with disjoint neighborhoods~$U$ we see that the codimension is unbounded and thus infinite.
\end{proof}

\section{Deformations of harmonic spinors}\label{sec.def.harm.spinors}

In the preceding section, we discussed the variation of the Dirac operator with respect to the Riemannian metric. We will now define the variation of the projection of the Dirac operator to the bundle spanned by small eigenspinors.

\subsection{The kernel projection}
Let $M$ be a closed connected spin manifold of dimension $m=2$ or $m=4$.
Let $g_0 \in \maR(M)$ have $\dim_\HH \ker\Dirac^{g_0}= |\alpha(M)|+ 2d$ for some $d\in \NN_0$, thus $\dim_\CC \ker\Dirac^{g_0}= 2|\alpha(M)|+ 4d$.
We choose $\epsilon>0$ such that \smash{$\spec\bigl(\Dirac^{g_0}\bigr)\cap(-\epsilon,\epsilon)\subset \{0\}$}.

Further, we consider an open star-shaped $C^1$-neighborhood~$\maU$ of~$g_0\in \maR_d(M)$ in the space of Riemannian metrics with respect to the $C^1$-topology,
such that \smash{$\epsilon \notin \spec\bigl(\Dirac^{g}\bigr)$} for all $g\in \maU$.
We take $\maU$ small enough, so that for all $g\in \maU$ the Dirac operator $\Dirac^g$ has precisely $2|\alpha(M)|+4d$ complex eigenvalues (counted with multiplicity) in $(-\epsilon,\epsilon)$. For $g\in \maU$, let
\[
K_g\ceq \bigoplus_{\lambda\in(-\epsilon,\epsilon)}\ker \bigl(\Dirac^g-\lambda\bigr) .
\]
Then $K\ceq \bigcup_{g\in \maU}K_g$ is a quaternionic vector bundle over $\maU$ of quaternionic rank $|\alpha(M)|+2d$, see, \eg \cite[Theorem~4.5.2]{nowaczyk:phd}.
The bundle $K \to \maU$ is a subbundle of finite rank of the bundle $L|_{\maU} \to \maU$, introduced in the previous section, and the $L^2$-orthonormal projection $L \to K$ is given by the spectral projection \smash{$\pi_{(-\epsilon,\epsilon)}^g$} to the sum of eigenspinors of $\Dirac^g$ with eigenvalues in $(-\epsilon,\epsilon)$.
The $L^2$-norm on $L$ gives a fiberwise real scalar product $\llangle\biargu\rrangle$ on $K$ that is invariant under multiplication by $a\in \HH$ with $|a|=1$, and it extends to a quaternionic scalar product $\llangle\biargu\rrangle_\HH$.

Be aware that $K$ is not a parallel subbundle of $L$, but -- as usual -- we obtain an induced connection on the bundle~$K$ by projecting to~$K$, i.e., \smash{$\pi_{(-\epsilon,\epsilon)}^g \circ \hat\nabla$} is the induced connection on~$K$. Then the inclusion $\iota_K\colon K\to L$ is not parallel, its derivative at $g_0$, applied to $\phi\in K$ is given~by
\[\biggl(\frac{\hat\nabla}{{\rm d}g}\bigg|_{g_0}\iota_K\biggr)\; \phi= \frac{\hat\nabla}{{\rm d}g}\bigg|_{g_0}\phi \;\;- \pi_{(-\epsilon,\epsilon)}^{g_0} \circ\frac{\hat\nabla}{{\rm d}g}\bigg|_{g_0}\phi= \pi_{(-\infty,-\epsilon]\cup [\epsilon,\infty)}^{g_0} \circ\frac{\hat\nabla}{{\rm d}g}\bigg|_{g_0} \phi ,\]
thus, the covariant derivative of $\iota_K$ at $g_0$ is a bilinear map from $\Gamma \bigl(\Sym T^*M\bigr) \times K$ to the orthogonal complement of $K$ in $L$.

The splitting $\Sigma^g M=\Sigma^g_+M\oplus \Sigma^g_-M$ and the metric $g\in \maU$ induce a splitting $K=K_+\oplus K_-$. These summands are quaternionic subbundles for $m=4$, but only complex subbundles for $m=2$. In the case $m=2$, multiplication by~$j$ maps $K_\pm$ to~$K_\mp$.

We now define the bundle $\maE \to \maU$ by
\[
\maE_g \ceq \End_{\omega\text{-ac}}^{\text{SA}} (K_g) ,
\]
that is, the fiber $\maE_g$ at $g$ is the space of self-adjoint endomorphisms of $K_g$ that anti-commute with the real volume element $\omega\ceq e_1\cliff e_2$ or $\omega\ceq e_1\cliff e_2\cliff e_3 \cliff e_4$ respectively, where $(e_1,\ldots,e_m)$ is a locally defined positively oriented orthonormal frame.
Then~$\maE$ inherits a connection, also denoted by $\hat\nabla$, from the connection on $K$, defined above.
When we restrict the Dirac operator~$\Dirac^g$ for each $g\in \maU$ to $K_g$, we obtain an endomorphism of $K_g$, and due to its self-adjointness and its anti-commutativity with the volume element, restricting to $K_g$ for each $g\in \maU$ yields a section of~$\maE$. Thus, for the inclusions $\iota_K\colon K\to L$ and $\iota_K^g\colon K_g\to L_g$
\[
\Dirac^{ \bullette}\big|_{K}=\Dirac^{ \bullette}\circ\iota_K\colon\ \maU \to \maE,
\qquad
g\mapsto \Dirac^g\big|_{K_g} = \Dirac^g\circ \iota_K^g,
\]
is a section of $\maE$. The zero locus of this section is $\maU\cap \maR_d(M)$. We denote the covariant linearization at $g_0$ of this map by
\begin{equation}\label{def.P}
P \ceq \pi_{(-\epsilon,\epsilon)}^{g_0} \circ\frac{\hat\nabla}{{\rm d}g}\bigg|_{g_0} \bigl(\Dirac^g \circ\iota_K^g\bigr)
\colon\ \Gamma \bigl(\Sym T^*M\bigr) \to \maE_{g_0} .
\end{equation}

In order to study the regularity and codimension of the zero level set $\maU\cap \maR_d(M)$ of
$\Dirac^{ \bullette}$, 
we will show that
$\rank \Dirac^{ \bullette} \ceq \dim_\RR \image P \leq 1$ 
at $g_0$ implies that the metric $g_0$ is conformal to a flat torus, so otherwise
$\rank \Dirac^{ \bullette} \geq 2$ 
which implies that $\codim (\maU\cap \maR_d(M)) \geq 2$, unless $g_0\in \Rcf(M)$.
This implies that the codimension of $\maR_d(M)$ in $\maR(M)$ is at least $2$.

\subsection{The energy-momentum tensor}

We define $\enerimp_{\phi,\psi}\in \Gamma\bigl(\Sym T^*M \otimes \HH\bigr)$ by the formula
\[
\enerimp_{\phi,\psi}(X,Y) \ceq
-\frac18(
\< X \cliff\na_Y \phi, \psi \>_\HH + \<Y \cliff \na_X \phi, \psi \>_\HH
+ \< \phi, X \cliff \na_Y \psi \>_\HH + \< \phi, Y \cliff\na_X \psi \>_\HH).
\]
For $p\in M$, $X,Y\in T_pM$, $\phi,\psi\in \ker\Dirac$, and a quaternion $\lambda$, we have
\[
\enerimp_{\phi\cdot \lambda,\psi}(X,Y) = \bar{\lambda}\cdot \enerimp_{\phi,\psi}(X,Y)
\text{ and }
\enerimp_{\phi,\psi\cdot \lambda}(X,Y)=\enerimp_{\phi,\psi}(X,Y)\cdot \lambda.
\]

The tensor $\enerimp_{\phi,\phi}$ is obtained by differentiating \smash{$\int_M\bigl\<\Dirac^g\phi,\phi\bigr\rangleR$} with respect to the Riemannian metric, see \cite[Section~6]{baer.gauduchon.moroianu:05} for details.
Such functionals often arise in action functionals in physics, and the differential is thus called an \define{energy-momentum tensor}.
For instance, it provides an energy-momentum term in Einstein's field equations arising from fermions.
We extend this terminology to the unsymmetrized variant $\enerimp_{\phi,\psi}$ used in this article.

\begin{Lemma}\label{lemma.em.TT}
If $\phi$ and $\psi$ are harmonic spinors, then $\tr \enerimp_{\phi,\psi}=0$.
\end{Lemma}

\begin{proof}
From the definition of $Q_{\phi,\psi}$, we have
\[
\enerimp_{\phi,\psi} (X,X)=
-\frac14 ( \< X \cliff \na_X \phi, \psi \>_\HH + \< \phi, X \cliff \na_X \psi \>_\HH ) .
\]
For an orthonormal basis $(e_1,\ldots,e_m)$ of the tangent space, we calculate
\[
\tr \enerimp_{\phi,\psi}
= \sum_{a=1}^m \enerimp_{\phi,\psi}(e_a,e_a)
= -\frac14 ( \< \Dirac\phi, \psi \>_\HH + \< \phi, \Dirac\psi \>_\HH )
= 0. \qedhere
\tag*{\qed}
\]
\renewcommand{\qed}{}
\end{proof}

\subsection{Variational formula}

We will now use the variational formula from Section~\ref{subsec:variation.of.D} in order to express the operator $P$ defined in \eqref{def.P} in terms of the energy-momentum tensor $\enerimp_{\phi,\psi}$.

\begin{Lemma}\label{lem.var.formula}
For $\phi,\psi \in \ker \Dirac^{g_0}$ and $h\in \Gamma\bigl(\Sym T^*M\bigr)$, we have
\[ \begin{split}
\llangle P(h)\phi,\psi\rrangleH
&=-\frac18 \int_M \sum_{a,b=1}^m(
\< e_a \cliff\na_{e_b} \phi, \psi \>_\HH + \< e_b \cliff\na_{e_a} \phi, \psi \>_\HH\\
&\hphantom{=-\frac18 \int_M \sum_{a,b=1}^m(}{}
+ \< \phi, e_a \cliff \na_{e_b} \psi \>_\HH + \< \phi, e_b\cliff \na_{e_a} \psi \>_\HH
)\cdot h(e_a,e_b) \dV^{g_0} \\
&= \llangle\enerimp_{\phi,\psi},h\rrangle.
\end{split}\]
\end{Lemma}

\begin{proof}
We compute
\begin{align*}
P
& = \pi_{(-\epsilon,\epsilon)}^{g_0} \circ \frac{\hat\nabla}{{\rm d}g}\bigg|_{g_0} \bigl(\Dirac^g \circ \iota_K^g \bigr)\\
&= \pi_{(-\epsilon,\epsilon)}^{g_0} \circ \biggl(\frac{\hat\nabla}{{\rm d}g}\bigg|_{g_0} \Dirac^g\biggr)\circ\iota_K^{g_0}
+ \pi_{(-\epsilon,\epsilon)}^{g_0} \circ \Dirac^{g_0}\circ\frac{\hat\nabla}{{\rm d}g}\bigg|_{g_0}\iota_k^g\\
&= \pi_{(-\epsilon,\epsilon)}^{g_0} \circ\biggl(\frac{\hat\nabla}{{\rm d}g}\bigg|_{g_0} \Dirac^g\biggr)\circ\iota_K^{g_0} .
\end{align*}
In the last equality, we used the fact that \smash{$\frac{\hat\nabla}{{\rm d}g}\big|_{g_0}\iota_k^g$} has values in the orthogonal complement of~$K$ in~$L$ which is the image of \smash{$\pi_{(-\infty,-\epsilon]\cup [\epsilon,\infty)}^{g_0}$}.
Thus \smash{$\Dirac^g\circ \frac{\hat\nabla}{{\rm d}g}\big|_{g_0}\iota_k^g$} is also in the image of \smash{$\pi_{(-\infty,-\epsilon]\cup [\epsilon,\infty)}^{g_0}$}. After postcomposing this with \raisebox{-1.5pt}{\smash{$\pi_{(-\epsilon,\epsilon)}^{g_0}$}}, we obtain~$0$.

Since $\phi,\psi \in \ker \Dirac^{g_0}$, this gives
\[
\llangle P(h)\phi,\psi\rrangleH
=
\biggl\langle\!\!\!\biggl\langle \pi_{(-\epsilon,\epsilon)}\circ\frac{\hat\nabla}{{\rm d}g}\bigg|_{g_0} \Dirac^g(h)\phi,\psi\biggr\rangle\!\!\!\biggr\rangle_\HH
=
\biggl\langle\!\!\!\biggl\langle \frac{\hat\nabla}{{\rm d}g}\bigg|_{g_0} \Dirac^g(h) \phi,\psi\biggr\rangle\!\!\!\biggr\rangle_\HH .
\]
Since $\Dirac^g$ is self-adjoint, its variation is also self-adjoint, and \eqref{first-var-D} gives us
\[ \begin{split}
\llangle P(h)\phi,\psi\rrangleH
={}&
\frac{1}{2} \biggl(
\biggl\langle\!\!\!\biggl\langle \frac{\hat\nabla}{{\rm d}g}\bigg|_{g_0} \Dirac^g(h) \phi,\psi\biggr\rangle\!\!\!\biggr\rangle_\HH
+ \biggl\langle\!\!\!\biggl\langle \phi, \frac{\hat\nabla}{{\rm d}g}\bigg|_{g_0} \Dirac^g(h) \psi\biggr\rangle\!\!\!\biggr\rangle_\HH
\biggr) \\
={}&
{-}\frac{1}{4} \int_M \sum_{a,b=1}^m \< e_a \cdot \nabla_{e_b} \phi,\psi \>_\HH h(e_a,e_b) \dV^{g_0}
-\frac{1}{8} \int_M \< \div^{g_0} h \cdot \phi,\psi \>_\HH \dV^{g_0} \\
&
{-}\frac{1}{4} \int_M \sum_{a,b=1}^m \< \phi,e_a \cdot \nabla_{e_b} \psi \>_\HH h(e_a,e_b) \dV^{g_0}
-\frac{1}{8} \int_M \< \phi,\div^{g_0} h \cdot \psi \>_\HH \dV^{g_0} \\
={}&
{-}\frac{1}{4} \int_M \sum_{a,b=1}^m ( \< e_a \cdot \nabla_{e_b} \phi,\psi \>_\HH
+ \< \phi, e_a \cdot \nabla_{e_b}\psi \>_\HH ) h(e_a,e_b) \dV^{g_0} .
\end{split} \]
By symmetrizing in $a$ and $b$, we conclude the lemma.
\end{proof}

\subsection{Calculations with the energy-momentum tensor}

\begin{Lemma}\label{lem.phipsi.quot.deriv}
Let $Z$ be a $($locally defined$)$ vector field with $|Z| = 1$. Then
\begin{align*}
|\psi|^2 \partial_Z |\phi|^2 - |\phi|^2 \partial_Z |\psi|^2
=
8\Re (\< Z \cliff \psi , \phi \>_\HH \cdot\enerimp_{\phi,\psi}(Z,Z)) .
\end{align*}
\end{Lemma}

\begin{proof}
Using \eqref{formula.HH.basis}, we compute
\begin{equation*}
|\psi|^2 \phi
= |Z \cliff \psi|^2 \phi
= Z \cliff \psi \< Z \cliff \psi ,\phi \>_\HH,
\end{equation*}
and with this
\begin{align*}
|\psi|^2 \partial_Z |\phi|^2
&= 2 \bigl\< |\psi|^2 \phi, \na_Z \phi \bigr\rangle_\RR = 2 \Re\bigl(\< |\psi|^2 \phi, \na_Z \phi \>_\HH\bigr) \\
&= 2 \Re(\< Z \cliff \psi \< Z \cliff \psi ,\phi \>_\HH , \na_Z \phi \>_\HH) \\
&= 2 \Re \bigl(\overline{\< Z \cliff \psi , \phi \>_\HH} \cdot\< Z \cliff \psi, \na_Z \phi \>_\HH \bigr) \\
&= -2 \Re \bigl(\overline{\< Z \cliff \psi , \phi \>_\HH} \cdot\< \psi, Z \cliff\na_Z \phi \>_\HH \bigr) \\
&= -2 \Re \bigl(\< Z \cliff \psi , \phi \>_\HH \cdot\overline{\< \psi, Z \cliff\na_Z \phi \>_\HH} \bigr) \\
&=
-2 \Re (
\< Z \cliff \psi,\phi \>_\HH \cdot \< Z \cliff \na_Z \phi, \psi \>_\HH
) .
\end{align*}
Similarly,
\begin{align*} 
|\phi|^2 \partial_Z |\psi|^2
&= -2 \Re \bigl(\overline{\< Z \cliff \phi,\psi \>_\HH}\cdot \< \phi, Z \cliff \na_Z \psi \>_\HH \bigr) \\
&= 2 \Re(\< Z \cliff \psi , \phi \>_\HH \cdot\< \phi, Z \cliff \na_Z \psi \>_\HH ).
\end{align*}
We finally get
\begin{align*} 
 |\phi|^2 \partial_Z |\psi|^2 - |\psi|^2 \partial_Z |\phi|^2
&
=
2 \Re (\< Z \cliff \psi , \phi \>_\HH\cdot
( \< Z \cliff \na_Z \phi, \psi \rangleH + \< \phi, Z \cliff \na_Z \psi \rangleH ))\\
&
= -8 \Re (\< Z \cliff \psi , \phi \rangleH\cdot
Q_{\phi,\psi}(Z,Z) ) .
\tag*{\qed}
\end{align*}
\renewcommand{\qed}{}
\end{proof}

Now we fix \smash{$\Phi\in \ker \Dirac^{g_0,+}\!\setminus\!\{0\}$} and \smash{$\Psi\in \ker \Dirac^{g_0,-}\!\setminus\!\{0\}$}.
On the set \smash{$U_\Psi := \{ x \in M \mid \Psi(x) \neq 0 \}$}, we set \smash{$G \ceq \grad \frac{|\Phi|^2}{|\Psi|^2}$}. Using Lemma~\ref{lem.phipsi.quot.deriv}, we compute
\begin{align}
g(G,Z)
=\partial_Z\biggl(\frac{|\Phi|^2}{|\Psi|^2}\biggr)
&=
\frac{|\Psi|^2 \partial_Z |\Phi|^2 - |\Phi|^2 \partial_Z |\Psi|^2}{|\Psi|^4} \nonumber\\
&=
\frac{8}{|\Psi|^4}
\Re (\< Z \cliff \Psi , \Phi \rangleH \cdot\enerimp_{\Phi,\Psi}(Z,Z)) .\label{eq.GZ-prod}
\end{align}

\section{Proofs of Theorems~\ref{thm.main.2} and~\ref{thm.main.3}}
\label{sec.proof.thm.main.2}

The proof of Theorem~\ref{thm.main.3} essentially relies on the following proposition.

\begin{Proposition}\label{prop.rank.two}
Let $(M,g_0)$ be a closed connected Riemannian spin manifold of dimension $m=2$ or $m=4$.
For any nonzero \smash{$\Phi\in \ker \Dirac^{g_0,+}$} and nonzero \smash{$\Psi\in \ker \Dirac^{g_0,-}$} with $\llangle\Phi,\Psi\rrangle_\HH=0$, the map
\begin{equation}\label{eq.quasi.Q}
 P_{\Phi,\Psi}\colon\ \Gamma\bigl(\Sym T^*M\bigr)\to \HH, \qquad h\mapsto \llangle P(h)\Phi,\Psi\rrangle_\HH
\end{equation}
has (real) rank at least $2$, unless $m=4$ and $(M,g_0)$ is conformal to a flat torus, \ie unless $g_0\in \Rcf(M)$.
\end{Proposition}

\begin{proof} 
We assume that the map $P_{\Phi,\Psi}$, defined in \eqref{eq.quasi.Q}, has rank at most $1$.
Due to Lemma~\ref{lem.var.formula}, this means that the image of
\[
\Gamma\bigl(\Sym T^*M\bigr)\to \HH, \qquad h\mapsto \llangle h,\enerimp_{\Phi,\Psi}\rrangle
= \int_M\<h,\enerimp_{\Phi,\Psi}\> \dvol^{g_0}
\]
is contained in $\RR \lambda$ for some $\lambda\in \HH$, $|\lambda|=1$. By replacing $\Psi$ by $\Psi \lambda^{-1}$, we obtain $\llangle h,\enerimp_{\Phi,\Psi}\rrangle\in \RR$ and all $h$. Thus $\enerimp_{\Phi,\Psi}$ is a (real-valued) section of $\Sym T^*M$.

From \eqref{eq.GZ-prod}, we obtain on $U_\Psi$
\begin{align*}
g(G, Z)
&=
\frac{8}{|\Psi|^4} \Re (\< Z \cliff \Psi , \Phi \rangleH\cdot \enerimp_{\Phi,\Psi}(Z,Z))\\
&=
\frac{8}{|\Psi|^4} \< Z \cliff \Psi , \Phi \rangleR \cdot \enerimp_{\Phi,\Psi}(Z,Z)\\
&=
\frac{8}{|\Psi|^4} \cdot \enerimp_{\Phi,\Psi}(Z,Z)\cdot g(V,Z) ,
\end{align*}
where $V\in \Gamma(TM)$ was defined as the vector field satisfying
\begin{equation*}
g(V,X) = \<X\cliff \Psi,\Phi\rangleR\qquad \text{for all } X\in TM.
\end{equation*}
This yields
\begin{equation} \label{eq.G-Fi}
g(G, Z) \cdot |Z|^2 = \frac{8}{|\Psi|^4} \enerimp_{\Phi,\Psi}(Z,Z) \cdot g(V,Z)
\end{equation}
for arbitrary $Z\in \Gamma(TM)$.

By choosing $Z\perp V$, we conclude that $G$ is of the form $fV$ for some $f\in C^\infty(M)$.
From \eqref{eq.G-Fi} we then obtain
\begin{equation*}
f\cdot g(V,Z) \cdot |Z|^2 = g(G, Z) \cdot |Z|^2 =\frac{8}{|\Psi|^4} \enerimp_{\Phi,\Psi}(Z,Z) \cdot g(V,Z)
\end{equation*}
for arbitrary $Z$, so for all $Z$ not perpendicular to $V$ we get
\begin{equation}\label{eq.q.G-Fi.shortened}
f \cdot g(Z,Z) = \frac{8}{|\Psi|^4} \enerimp_{\Phi,\Psi}(Z,Z).
\end{equation}
By a density argument, we see that \eqref{eq.q.G-Fi.shortened} holds for all $Z\in TM$. Taking the trace of both sides and using Lemma~\ref{lemma.em.TT}, we get{\samepage
\[
mf = f \cdot \tr g = \frac{8}{|\Psi|^4} \tr \enerimp_{\Phi,\Psi}=0,
\]
thus $\enerimp_{\Phi,\Psi} = 0$ and $G=fV=0$.}

In the case $m=4$, we conclude from $G=0$ that there is a constant $\rho>0$ so that $|\Phi|=\rho|\Psi|$ on~$U_\Psi$. We rescale $\Psi$ so that $|\Phi|=|\Psi|$ on~$U_\Psi$. Since we have $\enerimp_{\Phi,\Psi} = 0$, the sequence of Lem\-mas~8.2--8.7 in~\cite{maier:97}, followed by Lemmas~8.10--8.11 and Proposition~8.12 in \cite{maier:97}, tells us that~$(M,g_0)$ is conformal to a flat torus.

In the case $m=2$, we also conclude that there is a constant $\rho>0$ so that $|\Phi|=\rho|\Psi|$ on $U_\Psi$. Recall that on surfaces we can identify \smash{$\Sigma_-M\otimes_\C \Sigma_-M= \bigwedge^{1,0}M= (T_\C M)^*$}. It is thus a bundle with a holomorphic structure. Furthermore, multiplication with $j$ provides an isomorphism \smash{$ \overline{\Sigma_-M}\xrightarrow{\cong}\Sigma_+M$}. A spinor $\phi=\phi_++\phi_-$ with $\phi_\pm\in \Gamma(\Sigma_\pm M)$ is harmonic if and only if $\phi_+ \cdot j $ and $\phi_-$ are holomorphic sections of $\Sigma_-M$. We conclude that $M\setminus U_\Psi$ is finite. The condition $|\Phi|=\rho|\Psi|$ implies that there is a function $f\colon U_\Psi\to S^1\subset \CC$ such that
\[
\Phi\cdot j = \Psi\cdot \rho f
\]
on $U_\psi$. The holomorphicity of $\Phi\cdot j $ and $\Psi$ implies that $f$ is also holomorphic, the singularities in $M\setminus U_\Psi$ are removable, and finally we see that $f$ is constant, say $f=\lambda$. Thus $\Phi\cdot j = \Psi\cdot\rho \lambda$ and thus $\llangle\Phi,\Psi\rrangle_\HH = \|\Psi\|^2\cdot j\rho \bar\lambda\neq 0$ in contradiction to the assumption in the proposition.
\end{proof}

\begin{proof}[Proof of Theorem~\ref{thm.main.3}]
According to the preceding proposition, the map given in \eqref{eq.quasi.Q} has rank at least $2$. We may thus choose $S\subset \HH$, $\dim_\RR S=2$ in the image of this map.
Let $\pi^S$ be the orthogonal projection $\HH\to S$.

Let again $B_\epsilon(0)$ be the open \ball{} of radius $\epsilon$ in $\RR^2$. We extend $g_0$ to a map $g_{\bullette}\colon B_\epsilon(0)\to \maR(M)$, such that its differential at $0\in B_\epsilon(0)$ intersects with the kernel of the map $\pi^S\circ P_{\Phi,\Psi}$ given in~\eqref{eq.quasi.Q} only in $0$, \ie we have $\image \upd g_{\bullette}|_{0}\cap \ker \bigl(\pi^S\circ P_{\Phi,\Psi}\bigr)=\{0\}$.

We extend $\Phi$ and $\Psi$ to smooth sections of $K_+$ and $K_-$ along $g_\bullette$, \ie for all $(s,t)\in B_\epsilon(0)$ we have $\Phi_{s,t}\in K_+|_{g_{s,t}}$ and $\Psi_{s,t}\in K_-|_{g_{s,t}} $, depending smoothly on $(s,t)$. Theorem~\ref{thm.main.3} then follows from Proposition~\ref{prop.rank.two}.
\end{proof}

\begin{Remark}\label{rem.regularity}
The above proofs show Proposition~\ref{prop.rank.two} and Theorem~\ref{thm.main.3} for any regularity $C^\ell$ on $\Gamma\bigl(\Sym T^*M\bigr)$ and $\maR(M)$, $\ell\in \NN\cup\{\infty\}$, and the associated maps $P_{\Phi,\Psi}$, $g_\bullette$, $\Phi_\bullette$, $\Psi_\bullette$, etc.\ are continuously differentiable for any $\ell$.
\end{Remark}

\begin{proof}[Proof of Theorem~\ref{thm.main.2}]
Let $g_0\in \maR_d(M)$ for some $d\in \NN$, $d>0$.
In particular, in the case $m=4$, $(M,g_0)$ is not conformal to a flat torus.

In the case $m=4$, the condition $d>0$ implies that we have \smash{$\ker \Dirac^{g_0}_+\neq 0$} and \smash{$\ker \Dirac^{g_0}_-\neq 0$}. We choose arbitrary non-zero \smash{$\Phi\in \ker \Dirac^{g_0}_+$} and \smash{$\Psi\in \ker \Dirac^{g_0}_-$}. They satisfy \smash{$\llangle\Phi,\Psi\rrangle_\HH=0$} as \smash{$\ker \Dirac^{g_0}_+$} and \smash{$\ker \Dirac^{g_0}_-$} are orthogonal quaternionic subspaces. Thus $\Phi$ and $\Psi$ satisfy the assumptions of Theorem~\ref{thm.main.3}.

In the case $m=2$, $\Sigma_+M\to M$ is a complex line bundle, and thus the restriction of $\llangle\biargu\rrangle_\HH$ to $\Sigma_+M\to M$ takes values only in $\CC$, which yields a hermitian scalar product $\llangle\biargu\rrangle_\CC$. The condition $d>0$ implies that \smash{$\dim_\HH\ker \Dirac^{g_0}=\dim_\CC\ker \Dirac^{g_0}_+\geq 2$}. Thus we choose non-zero \smash{$\Phi,\widetilde\Phi\in \ker \Dirac^{g_0}_+$} with \smash{$0=\bigl\langle\!\!\bigl\langle \Phi,\widetilde\Phi\bigr\rangle\!\!\bigr\rangle_\CC=\bigl\langle\!\!\bigl\langle \Phi,\widetilde\Phi\bigr\rangle\!\!\bigr\rangle_\HH$}. Then $\Phi$ and \smash{$\Psi\ceq \widetilde\Phi\cdot j $} satisfy the assumptions of Theorem~\ref{thm.main.3}.

Assume $\ell<\infty$. Theorem~\ref{thm.main.3} allows us to apply the submersion theorem in the Banach manifold setting, see \cite[Theorem~3.5.4]{abraham.marsden.ratiu:1988} and conclude that $g_0$ has an open neighborhood $\maU$ in $\maR_{\leq d}(M)$,
such that
\begin{align*}
\maU\to S,\qquad g\mapsto
\pi^S\bigl(\llangle \Dirac^g \Phi, \Psi \rrangle_\HH\bigr)
\end{align*}
is a submersion.

This submersion maps $\maU\cap(\maR_d(M))$ to $0$. This implies that $\maR_d(M)$ is of codimension at least~$2$ in $\maR_{\leq d}(M)$ and thus in $\maR(M)$. The statement of Theorem~\ref{thm.main.2} follows.

For regularity $C^\infty$, we observe that the constructions lead to continuously differentiable maps
$g_\bullette \colon B_\epsilon(0)\to \maR^{C^\infty}(M)$
and to similar maps $\Phi_\bullette$ and $\Psi_\bullette$ mapping to spinors with $C^\infty$-regularity. Combined with the inclusion \smash{$\maR^{C^\infty}(M)\to\maR^{C^\ell}(M)$}, for $\ell<\infty$ and similar inclusions for spinors, the results follows for regularity $C^\infty$.
\end{proof}

\section{Proof of Theorem~\ref{thm.main}}\label{sec.proof.thm.main}

We define $\delta(g)$ as
\[
\delta(g)\ceq
\begin{cases}
 d & \text{if}\ g\in \maR_d(M),\\
\infty& \text{if}\ g\in \Rcf(M) .
\end{cases}
\]
This is an upper semi-continuous function on $\maR(M)$, \ie for a sequence $g_\alpha\to g$ we have $\delta(g)\geq \limsup_{\alpha\to \infty} \delta(g_\alpha)$. Theorem~\ref{thm.main} follows from Theorem~\ref{thm.main.2} by applying the following lemma to~$\maM\ceq \maR(M)$ with the function $\delta$ from above. Recall that $\Rcf(M)$ has infinite codimension in~$\maR(M)$ due to Proposition~\ref{prop.CF.codim}.

\begin{Lemma}
Let $\maM$ be an open and connected subset of a Banach space,
together with an upper semi-continuous function $\delta\colon \maM\to\overline\NN\ceq\NN\cup\{0,\infty\}$. For $d\in\ol\NN$, we write $\maM_d\ceq \delta^{-1}(d)$. We~assume that $\maM_0$ is open and dense in $\maM$ and that for $d>0$ every $\maM_d$ has codimension at least~$2$ in $\maM$ in the sense of the introduction. Then $\maM_0$ is connected.
\end{Lemma}

\begin{proof}
We define $\maM_{\leq d}\ceq \bigcup_{l\leq d}\maM_l=\delta^{-1}(\{1,2,\ldots,d\})$ and $\maM_{\leq \infty}\ceq \maM$.

Let us assume that $\maM_0$ is not connected.
Thus we find disjoint non-empty open subsets $\maU_1$ and $\maU_2$ of $\maM_0$ (and $\maM$) such that $\maU_1\cup \maU_2=\maM_0$.
As $\maM_0$ is dense, we have $\maM=\overline{\maM_0}=\overline{\maU_1}\cup \overline{\maU_2}$.
As $\maM$ is connected, it follows that $\maI\ceq \overline{\maU_1}\cap \overline{\maU_2}$ is non-empty.
Obviously $\maI\cap \maM_0=\varnothing$. Let $d_0\ceq\min\{\delta(\maI)\}\in \overline{\NN}$, and we assume $\delta(g_0)=d_0\geq 1$ for some $g_0\in \maI$.
By upper semi-continuity, there is an open neighborhood $\maV_0$ of $g_0$ in~$\maM$, such that $\maV_0\subset \maM_{\leq d_0}$.

As $\maM_{d_0}$ has codimension at least $2$ in $\maM$, there is a submanifold $\maN$ of codimension~$2$ of a~possibly smaller open neighborhood $\maV_1\subset \maV_0$ of $g$, such that $\maM_{d_0}\cap \maV_1\subset \maN$.

For any $\alpha\in \{1,2\}$, we see that $\maU_\alpha\cap \maV_1$ is open. It is non-empty as $g\in \overline{\maU_\alpha}$.
Thus $\maU_\alpha\cap \maV_1$ cannot be contained in the codimension $2$ submanifold $\maN$.
We choose $g_\alpha\in (\maU_\alpha\cap \maV_1)\setminus \maN$. Because $\maN$ has codimension $2$, we find a path \smash{$\gamma\colon[1,2]\to \maV_1\setminus \maN\subset \maM_{\leq (d_0-1)}$} with $\gamma(1)=g_1$ and $\gamma(2)=g_2$. Now \smash{$I_\alpha\ceq \gamma^{-1}\bigl(\ol{\maU_\alpha}\bigr)$} is a closed subset of $[1,2]$ containing $\alpha$ for $\alpha\in\{1,2\}$, and
\[
I_1\cup I_2= \gamma^{-1}\bigl(\ol{\maU_1}\cup \ol{\maU_2}\bigr)= \gamma^{-1}(\maM)=[1,2].
\]
Thus there is $t_0\in I_1\cap I_2$, \ie \smash{$\gamma(t_1)\in \ol{\maU_1}\cap \ol{\maU_2}=\maI$}. By construction of $d_0$, we have $\delta(\gamma(t_1))\geq d_0$ which contradicts $\gamma(t_1)\in\maM_{\leq (d_0-1)}$.
\end{proof}

\begin{proof}[Proof of Theorem~\ref{thm.main}]
If $\ell < \infty$, then Theorem~\ref{thm.main} follows from Theorem~\ref{thm.main.2} by applying the lemma to \smash{$\maM\ceq \maR^{C^\ell}(M)$} with the function $\delta$ from above.

In the case $\ell = \infty$, we have to face the problem that smooth functions constitute a Fr\'echet space and not a Banach space. One way to solve this is to replace the theorems on Banach spaces we used, by appropriate theorems on Fr\'echet spaces. More precisely, we could use versions of the submersion theorem, that still hold in the Fr\'echet space setting for $\ell=\infty$, see \cite[Proposition~III.11]{neeb.wagemann} and \cite[around Theorem~D]{gloeckner}.

However, it seems easier to us to argue as follows.
For a given pair of metrics in $\smash{\maR^{C^\infty}(M)}\subset\smash{\maR^{C^1}(M)}$, we already know
that there is a continuous path in \smash{$\maR^{C^1}(M)$} connecting this pair.
This continuous path can be approximated by a piecewise linear path of metrics, with break points consisting of $C^\infty$-metrics and sufficiently short linear pieces. Since the end points of this path and the break points are all $C^\infty$-metrics, the piecewise linear path is a path in $\maR^{C^\infty}(M)$, and Theorem~\ref{thm.main} follows in this case as well.
\end{proof}

\subsection*{Acknowledgements}
Bernd Ammann was supported by the CRC~1085 \emph{Higher Invariants} (Universit\"at Regensburg) and by SPP 2026 \emph{Geometry at infinity}, both funded by the DFG.

\pdfbookmark[1]{References}{ref}
\LastPageEnding

\end{document}